\documentclass{amsart}

\usepackage{amscd}
\usepackage{amssymb}
\usepackage{amsthm}
\usepackage{amsmath}

\usepackage{stmaryrd}

\newtheorem{theorem}{Theorem}[section]

\theoremstyle{plain}
\newtheorem{definition}[theorem]{Definition}

\newtheorem{prop}[theorem]{Proposition}
\newtheorem{cor}[theorem]{Corollary}
\newtheorem{lemma}[theorem]{Lemma}
\newtheorem{thm}[theorem]{Theorem}

\theoremstyle{definition}
\newtheorem*{stepone}{Case 1}
\newtheorem*{steptwo}{Case 2}
\newtheorem*{stepthree}{Case 3}
\newtheorem*{stepfour}{Case 4}
\newtheorem*{stepfive}{Case 5}

\theoremstyle{remark}

\newcommand{\OT}{\mathcal{O}_{T}}
\newcommand{\OC}{\mathcal{O}}

\newcommand{\N}{\mathbb{N}}
\newcommand{\Z}{\mathbb{Z}}

\newcommand{\C}{\mathbb{C}}
\newcommand{\K}{\mathbb{K}}

\newcommand{\g}{\mathfrak{g}}
\newcommand{\lb}{\mathfrak{b}}
\newcommand{\h}{\mathfrak{h}}
\newcommand{\p}{\mathfrak{p}}

\newcommand{\lam}{{\lambda}}

\newcommand{\Hom}{\mathrm{Hom}}
\newcommand{\End}{\mathrm{End}}
\newcommand{\Ext}{\mathrm{Ext}}

\numberwithin{equation}{section}

\begin{document}

\title{From Jantzen to Andersen Filtration via Tilting Equivalence}

\author{Johannes K\"ubel}
\address{Department of Mathematics, University of Erlangen, Germany}
\curraddr{Bismarckstr. 1 1/2, 91054 Erlangen, Germany}
\email{kuebel@mi.uni-erlangen.de}
\thanks{I would like to thank Wolfgang Soergel and Peter Fiebig for their support and many helpful discussions}


\date{December 5, 2010}


\begin{abstract}
The space of homomorphisms between a projective object and a Verma module in category $\OC$ inherits an induced filtration from the Jantzen filtration on the Verma module. On the other hand there is the Andersen filtration on the space of homomorphisms between a Verma module and a tilting module. Arkhipov's tilting functor, a contravariant self-equivalence of a certain subcategory of $\OC$, which maps projective to tilting modules induces an isomorphism of these kinds of Hom-spaces. We will show that this equivalence even identifies both filtrations.
\end{abstract}

\maketitle

\section{Introduction}

Let $\g \supset \lb \supset \h$ be a semisimple complex Lie algebra with a Borel and a Cartan. In BGG-category $\OC$, prominent objects are indecomposable projective and tilting modules. \\
Tilting modules were introduced in \cite{B} as selfdual Verma flag modules and the indecomposable tilting modules are classified by their highest weight. Let $\rho \in \h^{*}$ be the halfsum of positive roots relative to $\lb$ and let $T$ denote the ring of regular functions on the line $\C\rho$. $T$ is a quotient of the universal enveloping algebra of $\h$. For every weight $\lam \in \h^{*}$ the quotient map $\lb \twoheadrightarrow \h$ induces a $\lb$-module structure on $\C$. We denote this $\lb$-module by $\C_\lam$. The $\h$-module structure on $T$ also restricts to a $\lb$-module structure by the map $\lb \twoheadrightarrow \h$. Now we can form the Verma module $\Delta(\lam) = U(\g)\otimes_{U(\lb)} \C_\lam \in \g$-mod and the deformed Verma module
$$\Delta_{T}(\lam) = U(\g) \otimes_{U(\lb)} (\C_\lam \otimes T) \in \g-\mathrm{mod}-T$$
where tensor products without any specification are to be understood over $\C$. The $T$-module structure on $\Delta_T (\lam)$ is just multiplication from the right while the $\lb$-module structure on $\C_\lam \otimes T$ is the tensor product representation.\\
Taking the direct sum of the $T$-dual weight spaces and twisting the contragredient $\g$-module structure with a Chevalley automorphism leads to the deformed dual Verma module $\nabla_T(\lam) \in \g$-$\mathrm{mod}$-$T$. We will see that every invertible $T$-module homomorphism on the $\lam$-weight spaces extends to an injective homomorphism of $\g$-$T$-bimodules
$$\mathrm{can}:\Delta_T(\lam) \hookrightarrow \nabla_T(\lam)$$
which forms a basis of the $T$-module $\Hom_{\g-T}(\Delta_T(\lam),\nabla_T(\lam))$. Since $T$ can be understood as a polynomial ring in one variable $v$ we get the Jantzen filtration on $\Delta(\lam)$ by taking the images of $\mathrm{can}^{-1}(\nabla_T(\lam)v^{i})$ for $i=0,1,2,3,...$ under the surjection $\Delta_T(\lam) \twoheadrightarrow \Delta(\lam)$ induced by $\cdot \otimes_T \C$. For a projective object $P \in \OC$ we get an induced filtration on $\Hom_{\g}(P,\Delta(\lam))$.\\
Let $\widehat{T}$ be the completion of $T$ at the maximal ideal of 0. So we can identify $\widehat{T}$ with the ring of formal power series $\C \llbracket v \rrbracket$ in one variable.  In this article, we will also introduce deformed tilting modules which are certain $\g$-$\widehat{T}$-bimodules corresponding to tilting modules in $\OC$ after specializing with $\cdot \otimes_{\widehat{T}} \C$. Now let $K$ be such a deformed tilting module and consider the composition pairing
$$\Hom(\Delta_{\widehat{T}}(\lam),K) \times \Hom(K,\nabla_{\widehat{T}}(\lam)) \rightarrow \Hom(\Delta_{\widehat{T}}(\lam), \nabla_{\widehat{T}}(\lam))\cong {\widehat{T}}$$
where all Hom-spaces are meant to be homomorphisms of $\g$-${\widehat{T}}$-bimodules. We will see that this is a nondegenerate pairing of free $\widehat{T}$-modules of finite rank and leads to an injection
$$\Hom(\Delta_{\widehat{T}}(\lam),K) \hookrightarrow (\Hom(K,\nabla_{\widehat{T}}(\lam)))^{*}$$
where $(\cdot)^{*}$ denotes the ${\widehat{T}}$-dual.\\
Taking the preimages of the ${\widehat{T}}$-submodules $(\Hom(K,\nabla_{\widehat{T}}(\lam)))^{*}\cdot v^{i}$ under this embedding and applying $\cdot \otimes_{\widehat{T}} \C$ to these preimages defines the Andersen filtration on $\Hom_\g(\Delta(\lam),K\otimes_{\widehat{T}} \C)$.\\
In \cite{H}, Soergel introduces the tilting functor $t$ which forms a contravariant self-equivalence of the category of modules with a Verma flag, i.e. a filtration with subquotients isomorphic to Verma modules. The functor $t$ takes projective modules to tilting modules and sends a Verma module $\Delta(\nu)$ to the Verma module $\Delta(-2\rho-\nu)$. So $t$ induces an isomorphism of vector spaces
$$\Hom_\g(P,\Delta(\lam)) \stackrel{\sim}{\longrightarrow} \Hom_\g(\Delta(-2\rho-\lam),t(P))$$
which we denote by $t$ as well.
In this paper we will prove that $t$ even identifies the filtration induced by the Jantzen filtration on the left side with the Andersen filtration on the right side.\\

In \cite{I}, Soergel uses a hard Lefschetz argument to prove that the Andersen filtration on $\Hom_{\g}(\Delta(\lam),K)$ for $K$ a tilting module coincides with the grading filtration induced from the graded version of $\OC$ as described in \cite{J}. Since this is very similar to the result in \cite{A} about the semisimplicity of the subquotients of the Jantzen filtration, the relation of both filtrations might give an alternative proof of this semisimplicity.


\section{Preliminaries}

This section contains some results about the deformed category $\OC$ of a semisimple complex Lie algebra $\g$ with Borel $\lb$ and Cartan $\h$, which one can also find in \cite{C} and \cite{I}. By $S$ we will denote the universal enveloping algebra of the Cartan $\h$ which is equal to the ring of polynomial functions $\C[\h^{*}]$.
Let $T$ be a commutative, associative, noetherian, unital, local $S$-algebra with structure morphism $\tau: S \rightarrow T$. We call $T$ a \textit{local deformation algebra}.\\
In this article we will mostly deal with the $S$-algebras $R=S_{(0)}$, the localisation of $S$ at the maximal ideal of $0\in \h^{*}$, localisations $R_{\p}$ of $R$ at a prime ideal $\p$ of height 1 or the residue fields of these rings $\K_{\p}=R_{\p}/R_{\p}\p$. To apply results of this section to both filtrations we will also be concerned with the power series ring $\C \llbracket v \rrbracket$ in one variable and the quotient field $Q$ of $S$. All these rings are local deformation algebras.

\subsection{Deformed category $\OC$}

Let $T$ be a local deformation algebra with structure morphism $\tau: S\rightarrow T$ and let $M \in \g$-mod-$T$. For $\lam \in \h^{*}$ we set
$$M_\lam = \{m \in M| hm=(\lam + \tau)(h)m \, \, \forall h \in \h\}$$ 
where $(\lam+\tau)(h)$ is meant to be an element of $T$. We call the $T$-submodule $M_\lam$ the deformed $\lam$-weight space of $M$.\\

We denote by $\OT$ the full subcategory of all bimodules $M \in \g$-mod-$T$ such that $M= \bigoplus\limits_{\lam \in \h^{*}} M_\lam$ and with the properties that for every $m \in M$ the $\lb$-$T$-bimodule generated by $m$ is finitely generated as a $T$-module and that $M$ is finitely generated as a $\g$-$T$-bimodule. For example, if we put $T=\C$, $\OT$ is just the usual BGG-category $\OC$.\\
For $\lam \in \h^{*}$ we define the \textit{deformed Verma module}
$$\Delta_T(\lam) = U(\g) \otimes_{U(\lb)} T_\lam$$
where $T_\lam$ denotes the $U(\lb)$-$T$-bimodule $T$ with $\lb$-structure given by the composition $U(\lb) \rightarrow S \stackrel{\lam + \tau}{\longrightarrow} T$.\\

As in \cite{I}, we now introduce a functor 
$$d=d_\tau : \g \otimes T \mathrm{-mod} \longrightarrow \g \otimes T \mathrm{-mod}$$
by letting $dM \subset \Hom_T(M,T)^{\sigma}$ be the sum of all deformed weight spaces in the space of homomorphisms of $T$-modules from $M$ to $T$ with its $\g$-action twisted by an involutive automorphism $\sigma:\g \rightarrow \g$ with $\sigma|_{\h} = -\mathrm{id}$. We now set $\nabla_T(\lam)=d\Delta_T(\lam)$ for $\lam \in \h^{*}$ and call it the \textit{deformed nabla module}. As in \cite{I}, one shows $d\nabla_T(\lam) \cong \Delta_T(\lam)$ and that tensoring with a finite dimensional representation $E$ of $\g$ commutes with $d$ up to the choice of an isomorphism $dE \cong E$.

\begin{prop}[\cite{I}, Proposition 2.12.]\label{Soe1}
\begin{enumerate}
\item For all $\lam$ the restriction to the deformed weight space of $\lam$ together with the two canonical identifications $\Delta_T(\lam)_\lam \stackrel{\sim}{\rightarrow} T$ and $\nabla_T(\lam)_\lam \stackrel{\sim}{\rightarrow} T$ induces an  isomorphism
$$\Hom_{\OT}(\Delta_T(\lam),\nabla_T(\lam)) \stackrel{\sim}{\longrightarrow}T$$

\item For $\lam \neq \mu$ in $\h^{*}$ we have $\Hom_{\OT}(\Delta_T(\lam),\nabla_T(\mu))=0$.

\item For all $\lam,\mu \in \h^{*}$ we have $\Ext^{1}_{\OT}(\Delta_T(\lam),\nabla_T(\mu))=0$.
\end{enumerate}
\end{prop}

\begin{cor}[\cite{I}, Corollary 2.13.]\label{Soe2}
Let $M,N \in \OT$. If $M$ has a $\Delta_T$-flag and $N$ a $\nabla_T$-flag, then the space of homomorphisms $\Hom_{\OT}(M,N)$ is a finitely generated free $T$-module and for any homomorphism $T\rightarrow T'$ of local deformation algebras the obvious map defines an isomorphism
$$\Hom_{\OT}(M,N)\otimes_{T} T' \stackrel{\sim}{\longrightarrow}   \Hom_{\OC_{T'}}(M\otimes_{T} T',N\otimes_{T} T')$$
\end{cor}

\begin{proof}
This follows from Proposition \ref{Soe1} by induction on the length of the $\Delta_T$- and $\nabla_T$-flag.
\end{proof}

If $\mathfrak{m}\subset T$ is the unique maximal ideal in our local deformation algebra $T$ we set $\K=T/\mathfrak{m}T$ for its residue field.

\begin{thm}[\cite{C}, Propositions 2.1 and 2.6] \label{Fie1}

\begin{enumerate}
\item The base change $\cdot \otimes_T \K$ gives a bijection\\
\begin{displaymath}
		\begin{array}{ccc}
		
			\left\{\begin{array}{c}
        \textrm{simple isomorphism}\\
      	\textrm{classes of $\OT$}
    	\end{array}\right\}
		&
		\longleftrightarrow
		&
			\left\{\begin{array}{c}
        \textrm{simple isomorphism}\\
      	\textrm{classes of $\OC_{\K}$}
    	\end{array}\right\}
		\end{array}
	\end{displaymath}

\item The base change $\cdot \otimes_T \K$ gives a bijection\\	
\begin{displaymath}
		\begin{array}{ccc}
		
			\left\{\begin{array}{c}
        \textrm{projective isomorphism}\\
      	\textrm{classes of $\OT$}
    	\end{array}\right\}
		&
		\longleftrightarrow
		&
			\left\{\begin{array}{c}
        \textrm{projective isomorphism}\\
      	\textrm{classes of $\OC_{\K}$}
    	\end{array}\right\}
		\end{array}
	\end{displaymath}
\end{enumerate}
\end{thm}

The category $\OC_{\K}$ is the direct summand of the category $\OC$ over the Lie algebra $\g \otimes \K$ consisting of all objects whose weights lie in the complex affine subspace $\tau + \h^{*} = \tau + \Hom_\C(\h,\C) \subset \Hom_\K(\h\otimes \K,\K)$ for $\tau$ the restriction to $\h$ of the map that makes $\K$ to a $S$-algebra. So the simple objects of $\OC_\K$ as well as the ones of $\OT$ are parametrized by their highest weight in $\h^{*}$. Denote by $L_T(\lam)$ the simple object with highest weight $\lam$. We also use the usual partial order on $\h^{*}$ to partially order $\tau + \h^{*}$.

\begin{thm}[\cite{C}, Propositions 2.4 and 2.7] \label{Fie2}
Let $T$ be a local deformation algebra and $\K$ its residue field. Let $L_T(\lam)$ be a simple object in $\OT$.

\begin{enumerate}
\item There is a projective cover $P_T(\lam)$ of $L_T(\lam)$ in $\OT$. Every projective object in $\OT$ is isomorphic to a direct sum of projective covers.
\item $P_T(\lam)$ has a Verma flag, i.e. a finite filtration with subquotients isomorphic to Verma modules, and for the multiplicities we have the BGG-reciprocity formula 
$$(P_T(\lam):\Delta_T(\mu)) = [\Delta_{\K}(\mu):L_{\K}(\lam)]$$
for all Verma modules $\Delta_T(\mu)$ in $\OT$.
\item Let $T \rightarrow T'$ be a homomorphism of local deformation algebras and $P$ projective in $\OT$. Then $P\otimes_T T'$ is projective in $\OC_{T'}$ and the natural transformation
$$\Hom_{\OT}(P,\cdot) \otimes_T T' \longrightarrow \Hom_{\OC_{T'}}(P \otimes_T T', \cdot \otimes_T T')$$
is an isomorphism of functors from $\OT$ to $T'$-mod.
\end{enumerate}
\end{thm}

Since Fiebig works over a complex symmetrizable Kac-Moody algebra, he has to introduce truncated subcategories and has to set some finiteness assumptions. In case of a finite dimensional semisimple Lie algebra we do not need these technical tools.

\subsection{Block decomposition}

Let $T$ again denote a local deformation algebra and $\K$ its residue field.

\begin{definition}
Let $\sim_T$ be the equivalence relation on $\h^{*}$ generated by $\lam \sim_T \mu$ if $[\Delta_{\K}(\lam):L_{\K}(\mu)] \neq 0$.
\end{definition}

\begin{definition}
Let $\Lambda \in \h^{*}/\sim_{T}$ be an equivalence class. Let $\OC_{T,\Lambda}$ be the full subcategory of $\OT$ consisting of all modules $M$ such that every highest weight of a subquotient of $M$ lies in $\Lambda$.
\end{definition}

\begin{prop}[\cite{C}, Proposition 2.8] \label{Fie3}
The functor
\begin{eqnarray*}
		\begin{array}{ccc}
		\bigoplus\limits_{\Lambda \in \h^{*}/\sim_T} \OC_{T,\Lambda}  &\longrightarrow&  \OT\\
		(M_\Lambda)_{\Lambda \in \h^{*}/\sim_T}& \longmapsto & \bigoplus\limits_{\Lambda \in \h^{*}/\sim_T} M_\Lambda
		
		\end{array}
	\end{eqnarray*}
is an equivalence of categories.
\end{prop}
The isomorphism above is called \textit{block decomposition}.\\
Later we will be especially interested in the case $T=R=S_{(0)}$ where $S_{(0)}$ denotes the localisation of $S$ at the maximal ideal generated by $\h$, i.e. the maximal ideal of $0 \in \h^{*}$.
Since $\sim_R = \sim_\C$, the block decomposition of $\OC_R$ corresponds to the block decomposition of the BGG-category $\OC$ over $\g$.\\

Let $\tau:S \rightarrow \K$ be the induced map that makes $\K$ into a $S$-algebra. Restricting to $\h$ and extending with $\K$ yields a $\K$-linear map $\h \otimes \K \rightarrow \K$ which we will also call $\tau$. Let $\mathcal{R} \supset \mathcal{R}^{+}$ be the root system with positive roots according to our data $\g\supset \lb \supset \h$. For $\lam \in \h_{\K}^{*}=\Hom_\K(\h \otimes \K,\K)$ and $\check{\alpha} \in \h$ the dual root of a root $\alpha \in \mathcal{R}$ we set $\left\langle \lam, \check{\alpha}\right\rangle_{\K} = \lam(\check{\alpha}) \in \K$. Let $\mathcal{W}$ be the Weyl group of $(\g,\h)$.

\begin{definition}
For $\mathcal{R}$ the root system of $\g$ and $\Lambda \in \h^{*}/\sim_T$ we define 
$$\mathcal{R}_T(\Lambda)=\{\alpha \in \mathcal{R} | \left\langle \lam + \tau, \check{\alpha}\right\rangle_{\K} \in \Z \subset \K \text{ for some } \lam \in \Lambda\}$$
and call it the integral roots corresponding to $\Lambda$. Let $\mathcal{R}_T^{+}(\Lambda)$ denote the positive roots in $\mathcal{R}_T(\Lambda)$ and set $$\mathcal{W}_T(\Lambda)=\langle \{s_\alpha \in \mathcal{W} | \alpha \in \mathcal{R}_T^{+}(\Lambda)\}\rangle \subset \mathcal{W}$$
We call it the integral Weyl group with respect to $\Lambda$.
\end{definition}
From \cite{C} Corollary 3.3 it follows that
$$\Lambda=\mathcal{W}_T(\Lambda)\cdot \lam \text{  for any  } \lam \in \Lambda$$
where we denote by $\cdot$ the $\rho$-shifted dot-action of the Weyl group.\\

Since most of our following constructions commute with base change, we are particularly interested in the case when $T=R_{\p}$ is a localization of $R$ at a prime ideal $\p$ of height one. Applying the functor $\cdot \otimes_R R_{\p}$ will split the deformed category $\OT$ into generic and subgeneric blocks which is content of the next

\begin{lemma}[\cite{D}, Lemma 3] \label{Fie4}
Let $\Lambda \in \h^{*}/\sim_{R}$ and let $\p \in R$ be a prime ideal.
\begin{enumerate}
\item If $\check{\alpha} \notin \p$ for all roots $\alpha \in \mathcal{R}_{R}(\Lambda)$, then $\Lambda$ splits under $\sim_{R_\p}$ into generic equivalence classes.
\item If $\p = R\check{\alpha}$ for a root $\alpha \in \mathcal{R}_{R}(\Lambda)$, then $\Lambda$ splits under $\sim_{R_\p}$ into subgeneric equivalence classes of the form $\{\lam,s_{\alpha}\cdot \lam\}$.
\end{enumerate}
\end{lemma}
We recall that we denote by $P_T(\lam)$ the projective cover of the simple object $L_T(\lam)$. It is indecomposable and up to isomorphism uniquely determined.
For an equivalence class $\Lambda \in \h^{*}/\sim_T$ which contains $\lam$ and is generic, i.e. $\Lambda=\{\lam\}$, we get $P_T(\lam)=\Delta_T(\lam)$. If $\Lambda=\{\lam,\mu\}$ and $\mu< \lam$, we have $P_T(\lam)=\Delta_T(\lam)$ and there is a non-split short exact sequence in $\OT$
$$0\rightarrow \Delta_T(\lam) \rightarrow P_T(\mu) \rightarrow \Delta_T(\mu)\rightarrow 0$$ 
In this case, every endomorphism $f: P_T(\mu) \rightarrow P_T(\mu)$ maps $\Delta_T(\lam)$ to $\Delta_T(\lam)$ since $\lam>\mu$. So $f$ induces a commutative diagram 
\begin{eqnarray*}
\begin{CD}
   0   @>>> \Delta_{T}(\lam) @>>> P_{T}(\mu) @>>> \Delta_{T}(\mu) @>>> 0\\
   @VVV @V f_\lam VV @VV f V @VV f_\mu V @VVV\\
   0   @>>> \Delta_{T}(\lam) @>>> P_{T}(\mu) @>>> \Delta_{T}(\mu) @>>> 0 
\end{CD}\end{eqnarray*}
Since endomorphisms of Verma modules correspond to elements of $T$, we get a map
\begin{displaymath}
		\begin{array}{ccc}
		\chi: \End_{\OC_{T}}(P_{T}(\mu))& \longrightarrow & T \oplus T \\
		f & \longmapsto &(f_\lam , f_\mu)  
		\end{array}
\end{displaymath}

For $\p=R\check{\alpha}$ we define $R_\alpha :=R_\p$ for the localization of $R$ at the prime ideal $\p$.

\begin{prop}[\cite{C}, Corollary 3.5] \label{Fie5}
Let $\Lambda \in \h^{*}/\sim_{R_\alpha}$.
If $\Lambda=\{\lam,\mu\}$ and $\lam=s_{\alpha}\cdot \mu >\mu$, the map $\chi$ from above induces an isomorphism of $R_\alpha$-modules
$$\End_{\OC_{R_\alpha}}(P_{R_\alpha}(\mu)) \cong \left\{(t_\lam,t_\mu) \in {R_\alpha}\oplus {R_\alpha} \middle| t_\lam \equiv t_\mu \text{ mod } \check{\alpha}\right\}$$
\end{prop}

\section{Tilting modules and tilting equivalence} 

In this chapter, $T$ will be a localisation of $R=S_{(0)}$ at a prime ideal $\p \subset R$ and let $\K$ be its residue field. Let $\lam \in \h^{*}$ be such that $\Delta_{\K}(\lam)$ is a simple object in $\OC_{\K}$. Thus, we have $\Delta_{\K}(\lam) \cong \nabla_{\K}(\lam)$ and the canonical inclusion $\mathrm{Can}:\Delta_T(\lam) \hookrightarrow \nabla_T(\lam)$ becomes an isomorphism after applying $\cdot \otimes_T \K$. So by Nakayama's lemma, we conclude that Can was bijective already.

\subsection{Deformed tilting modules}

\begin{definition}
By $\mathcal{K}_T$ we denote the full subcategory of $\OT$ which
\begin{enumerate}
\item includes the self-dual deformed Verma modules
\item is stable under tensoring with finite dimensional $\g$-modules
\item is stable under forming direct sums and summands. 
\end{enumerate}
\end{definition}
For $T=S/S\h=\C$ the category $\mathcal{K}_T$ is just the usual subcategory of tilting modules of the category $\OC$ over $\g$. In general, $\mathcal{K}_\K$ is the category of tilting modules of category $\OC$ over the Lie algebra $\g \otimes \K$ whose weights live in the affine complex subspace $\tau + \h^{*} \subset \Hom_\K(\h \otimes \K,\K)$, where $\tau:\h\otimes {\K} \rightarrow \K$ comes from the map that makes $\K$ into an $S$-algebra.

\begin{prop}\label{tilt1}
The base change $\cdot \otimes_T \K$ gives a bijection\\	
\begin{displaymath}
		\begin{array}{ccc}
		
			\left\{\begin{array}{c}
        \textrm{isomorphism classes}\\
      	\textrm{of $\mathcal{K}_T$}
    	\end{array}\right\}
		&
		\longleftrightarrow
		&
			\left\{\begin{array}{c}
        \textrm{isomorphism classes}\\
      	\textrm{of $\mathcal{K}_\K$}
    	\end{array}\right\}
		\end{array}
	\end{displaymath}
\end{prop}

\begin{proof}
For $K,H \in \mathcal{K}_T$ with $K \otimes_T \K \cong H\otimes_T \K$ we conclude $K \cong H$ from Nakayama's lemma applied to the weight spaces, since the weight spaces of tilting modules are finitely generated and free over $T$. This shows injectivity.\\
For surjectivity we only have to show that every indecomposable tilting module in $\OC_\K$ has an indecomposable preimage in $\mathcal{K}_T$. Since we are not working over a complete local ring we cannot apply the idempotent lifting lemma as in the proof of Proposition 3.4. in \cite{I}. Rather, the proof works very similar to the proof of Theorem 6 in \cite{F}.\\
Let $K \in \OC_\K$ be an indecomposable tilting module. For the sake of simplicity, we will assume the highest weight of $K$ to be regular. The singular case is treated analogously.\\
If the highest weight $\lam$ of $K$ is minimal in its equivalence class under $\sim_T$ then, $K\cong \Delta_\K(\lam) \cong \nabla_\K(\lam)$ and we can take $\Delta_T(\lam)$ as a preimage of $K$ in $\mathcal{K}_T$.\\
Now denote by $\overline{\lam}$ the equivalence class of $\lam$ in $\h^{*}/\sim_T$ and let $\lam = w \cdot \mu$ with $w \in \mathcal{W}_T(\overline{\lam})$ and $\mu$ minimal in $\mathcal{W}_T(\overline{\lam}) \cdot \lam$. In addition, let $w=s_1s_2...s_n$ be a minimal expression of $w$ with simple reflections $s_i \in \mathcal{W}_T(\bar{\lam})$. We denote by $\theta_i=\theta_i^{\K}$ the translation functor of $\OC_\K$ through the $s_i$-wall. Then $K$ is a direct summand of the tilting module 
$$M = \theta_1...\theta_n \Delta_\K(\mu)$$
and we get a decomposition $M \cong K \oplus K'$. $K'$ decomposes into indecomposable tilting modules with highest weights of the form $w'\cdot \mu$ and $l(w)>l(w')$ where $l$ denotes the length of a Weyl group element. By using induction on the length of $w$, we get a preimage $\tilde{K} \in \mathcal{K}_T$ of $K'$ together with a splitting inclusion
$$\tilde{K} \otimes_T \K \hookrightarrow M$$
Using Nakayama's lemma, this induces a splitting lift
$$\tilde{K} \hookrightarrow \theta_1^{T}...\theta_n^{T} \Delta_T(\mu)$$
where $\theta_i^{T}$ denotes the $T$-deformed wall-crossing functor of $\OT$ corresponding to $\theta_i^{\K}$.
Finally, the cokernel of this inclusion is the indecomposable tilting module in $\mathcal{K}_T$ we were looking for.
\end{proof}

\subsection{Tilting functor}

Let $\mathcal{M}_T$ denote the subcategory of all modules in $\OT$ admitting a Verma flag. We recall the map $\tau:S \rightarrow T$ which makes $T$ into a $S$-algebra. We get a new $S$-algebra structure on $T$ via $\tau \circ \gamma: S \rightarrow T$, where $\gamma:S \rightarrow S$ is the isomorphism given by $\gamma(h)=-h$ for all $h \in \h$. We denote this new $S$-algebra by $\overline{T}$. Let $S_{2\rho}$ be the semiregular $U(\g)$-bimodule of \cite{H}. If $N$ is a $\g \otimes T$-module which decomposes into weight spaces we get a $\overline{T}$-module 
$$N^{\star}= \bigoplus\limits_{\lam \in \h^{*}} \Hom_T(N_\lam,T)$$
Then we get a $\g \otimes \overline{T}$-module via the $\g$-action
$(Xf)(v)=-f(Xv)$
for all $X \in \g, f\in N^{\star}$ and $v\in N$.\\
Now we get a functor 
$$t'_T:\mathcal{M}_T \longrightarrow \mathcal{M}_{\overline{T}}^{opp}$$
by setting $t'_T(M) = (S_{2\rho} \otimes_{U(\g)}M)^{\star}$.\\
For $T$ a localisation of $S$ at a prime ideal $\p$ which is stable under $\gamma$, for a residue field of this or for the ring $\C \llbracket v \rrbracket$ of formal power series coming from a line $\C\lam \subset \h^{*}$, $\gamma$ induces an isomorphism of $S$-algebras $\gamma:T \stackrel{\sim}{\longrightarrow} \overline{T}$ which induces an equivalence of categories
$$\gamma:\mathcal{M}_T \longrightarrow \mathcal{M}_{\overline{T}}$$

\begin{thm}[\cite{D}, Section 2.6] \label{Fie6}
The functor $t_T=  (S_{2\rho} \otimes_{U(\g)} \cdot)^{\star} \circ \gamma$ induces an equivalence of categories
$$t_T : \mathcal{M}_T \longrightarrow \mathcal{M}_T^{opp}$$
which respects block decomposition, makes short exact sequences to short exact sequences and sends a Verma module $\Delta_T(\lam)$ to the Verma module $\Delta(-2\rho-\lam)$ for any weight $\lam \in \h^{*}$.
\end{thm}

\begin{prop}
Let $\lam \in \h^{*}$. Then
$$t_T(P_T(\lam)) \cong K_T(-2\rho -\lam)$$
where $P_T(\lam)$ denotes the indecomposable projective cover of $L_T(\lam)$ and $K_T(\mu)$ the up to isomorphism unique indecomposable deformed tilting module with highest weight $\mu \in \h^{*}$.
\end{prop}

\begin{proof}
We set $\mu = -2 \rho - \lam$. This proof is very similar to the proof of Proposition 3.1 of \cite{G}. As we have already seen, the tilting module $K_T(\mu)$ can be described as the up to isomorphism unique indecomposable module in $\OT$ with the properties:
\begin{enumerate}
\item $K_T(\mu)$ admits a Verma flag
\item $K_T(\mu)$ has a $\nabla_T$-flag
\item $K_T(\mu)_\mu$ is free of rank one over $T$
\item If $\gamma$ is a weight of $K_T(\mu)$, we have $\gamma \leq \mu$.
\end{enumerate}
Since $t_T$ is fully faithful, we already conclude the indecomposability of $t_T(P_T(\lam))$. Theorem \ref{Fie6} also tells us that $t_T(P_T(\lam))$ has a Verma flag.
By BGG-reciprocity we get a short exact sequence
$$N \hookrightarrow P_T(\lam) \twoheadrightarrow \Delta_T(\lam)$$
where $N$ has a Verma flag in which the occurring weights are strictly larger than $\lam$. By applying $t_T$ we get a new short exact sequence
$$t_T(N) \twoheadleftarrow t_T(P_T(\lam)) \hookleftarrow t_T(\Delta_T(\lam))$$
By induction on the length of the Verma flag of $N$ we conclude that the weights of $t_T(N)$ are strictly smaller than $\mu$. So the weight space $(t_T(P_T(\lam)))_\mu$ is free of rank one, since $t_T(\Delta_T(\lam))\cong \Delta_T(\mu)$. Now, $P_T(\lam)$ being projective and $t_T$ being fully faithful, we get
$$\Ext^{1}_{\OT}(\Delta_T(\delta), t_T(P_T(\lam)))= 0 \,  \,\,\, \forall \, \delta \in \h^{*}$$
Now we set $D_T=t_T(P_T(\lam))$ and consider the diagramm
\begin{eqnarray*}
	\begin{array}{ccccc}
	\Delta_T(\mu)&\hookrightarrow& K_T(\mu)& \twoheadrightarrow& \text{coker}\\
	\parallel&&&&\\
	\Delta_T(\mu)&\hookrightarrow& D_T& \twoheadrightarrow& \text{coker}'
	\end{array}
\end{eqnarray*}

Since coker has a Verma flag we conclude $\Ext^{1}_{\OT}(\text{coker},D_T)=0$ using induction on the length of a Verma flag of coker. So the restriction induces a surjection
$\Hom_{\OT}(K_T(\mu),D_T) \twoheadrightarrow \Hom_{\OT}(\Delta_T(\mu),D_T)$
and we get a map $\alpha: K_T(\mu) \rightarrow D_T$ which induces the identity on $\Delta_T(\mu)$. For the same reason we also get a map $\beta: D_T \rightarrow K_T(\mu)$ with the same property, since coker$'$ has a Verma flag while $K_T(\mu)$ admits a nabla flag which implies $\Ext_{\OT}^{1}(\text{coker}',K_T(\mu))=0$ by Propostion \ref{Soe1}.
Applying the base change functor $\cdot \otimes_T \K$, we get two maps\\  $\varphi:=(\beta \circ \alpha)\otimes \text{id}_{\K}:K_{\K}(\mu) \rightarrow K_{\K}(\mu)$ and $\psi:=(\alpha \circ \beta)\otimes \text{id}_{\K}:D_{T}\otimes_T \K \rightarrow D_{T}\otimes_T \K$ which induce the identity on $\Delta_\K(\mu)$. We conclude that $\varphi$ and $\psi$ are not nilpotent and since $D_T  \otimes_T \K$ and $K_\K(\mu)$ have finite length and are indecomposable, it follows from the Fitting lemma that $\varphi$ and $\psi$ are isomorphisms. Now by Nakayama's lemma applied to all weight spaces, we conclude $K_T(\mu) \cong D_T$.
\end{proof}

\begin{lemma}\label{tilt}
Let $T \rightarrow T'$ be a homomorphism of $S$-algebras where also $T'$ is a localisation of S at a prime ideal which is stable under $\gamma$, its residue field or the ring of formal power series $\C \llbracket v \rrbracket$. Let $M,N \in \mathcal{M}_T$ and let $M$ be projective in $\OT$. Then the diagram
\begin{eqnarray*}
\begin{CD}
   \Hom_{\OT}(M,N)\otimes_T T'   @>>> \Hom_{\mathcal{O}_{T'}}(M\otimes_T T',N\otimes_T T')\\
   @VVt_T \otimes \mathrm{id}_{T'} V @VVt_{T{'}} V\\
   \Hom_{\OT}(t_T(N),t_T(M))\otimes_T T' @>>> \Hom_{\mathcal{O}_{T'}}(t_{T'}(N\otimes_T T'),t_{T'}(M\otimes_T T'))
\end{CD}
\end{eqnarray*}
commutes, where the horizontals are the base change isomorphisms and the verticals are induced by the tilting functors $t_T$ resp. $t_{T'}$.
\end{lemma}

\begin{proof}
All composition factors of the tilting functor commute with base change in the sense of the lemma.
\end{proof}

\section{The Jantzen and Andersen filtrations}

We fix a deformed tilting module $K \in \mathcal{K}_T$ and let $\lam \in \h^{*}$. The composition of homomorphisms induces a $T$-bilinear pairing
\begin{eqnarray*}
	\begin{array}{ccc}
	\Hom_{\OT}(\Delta_T(\lam),K)\times \Hom_{\OT}(K,\nabla_T(\lam)) & \longrightarrow & \Hom_{\OT}(\Delta_T(\lam),\nabla_T(\lam))\cong T\\
	\\
	(\varphi,\psi) & \longmapsto & \psi \circ \varphi
	\end{array}
\end{eqnarray*}

For any $T$-module $H$ we denote by $H^{*}$ the $T$-module $\Hom_T(H,T)$. As in \cite{I} Section 4 one shows that for $T$ a localization of $S$ at a prime ideal $\p$ or for $T=\C \llbracket v \rrbracket$ our pairing is nondegenerate and induces an injective map
$$E=E_T^{\lam}(K): \Hom_{\OT}(\Delta_T(\lam),K) \longrightarrow \left(\Hom_{\OT}(K,\nabla_T(\lam))\right)^{*}$$
of finitely generated free $T$-modules.\\
If we take $T=\C \llbracket v \rrbracket$ the ring of formal power series around the origin on a line $\C\delta \subset \h^{*}$ not contained in any hyperplane corresponding to a reflection of the Weyl group, we get a filtration on $\Hom_{\OT}(\Delta_T(\lam),K)$ by taking the preimages of $\left (\Hom_{\OT}(K,\nabla_T(\lam)) \right )^{*}\cdot v^{i}$ for $i=0,1,2,...$ under $E$.

\begin{definition}[\cite{I}, Definition 4.2.]
Given $K_\C \in \mathcal{K}_\C$ a tilting module of $\OC$ and $K \in \mathcal{K}_{\C \llbracket v \rrbracket}$ a preimage of $K_\C$ under the functor $\cdot \otimes_{\C \llbracket v \rrbracket} \C$, which is possible by Proposition \ref{tilt1} with $S\rightarrow \C \llbracket v \rrbracket$ the restriction to a formal neighbourhood of the origin in the line $\C \rho$, then the image of the filtration defined above under specialization $\cdot \otimes_{\C \llbracket v \rrbracket} \C$ is called the \textit{Andersen filtration} on $\Hom_\g (\Delta(\lam), K_\C)$.
\end{definition}

The Jantzen filtration on a Verma module $\Delta(\lam)$ induces a filtration on the vector space $\Hom_\g(P , \Delta(\lam))$, where $P$ is a projective object in $\OC$. Now consider the embedding $\Delta_{\C \llbracket v \rrbracket}(\lam) \hookrightarrow \nabla_{\C \llbracket v \rrbracket}(\lam)$. Let $P_{\C \llbracket v \rrbracket}$ denote the up to isomorphism unique projective object in $\OC_{\C \llbracket v \rrbracket}$ that maps to $P$ under $\cdot \otimes_{\C \llbracket v \rrbracket}\C$, which is possible by Theorem \ref{Fie1}. Then we get the same filtration by taking the preimages of $\Hom_{\OC_{\C \llbracket v \rrbracket}}(P_{\C \llbracket v \rrbracket},\nabla_{\C \llbracket v \rrbracket}(\lam)) \cdot v^{i}$, $i=0,1,2,...$, under the induced inclusion
$$J=J_T^{\lam}(P):\Hom_{\OT}(P_T,\Delta_T(\lam)) \longrightarrow \Hom_{\OT}(P_T,\nabla_T(\lam))$$
for $T=\C \llbracket v \rrbracket$ and taking the images of these filtration layers under the map $\Hom_{\OT}(P_T,\Delta_T(\lam))\twoheadrightarrow \Hom_{\g}(P,\Delta(\lam))$ induced by $\cdot \otimes_T \C$. \\

For what follows, we define $\mu'=-2\rho-\mu$ and $\lam'=-2\rho-\lam$. To avoid ambiguity, we sometimes write $(\cdot)^{*_T}$, when we mean the $T$-dual of the $T$-module in brackets.

\begin{thm}
Let $\lam, \mu \in \h^{*}$. Denote by $R=S_{(0)}$ the localization of $S$ at $0$. There exists an isomorphism $L=L_R(\lam,\mu)$ which makes the diagram
\begin{eqnarray}\label{eqn:DiaS}
	\begin{CD}
   \Hom_{\OC_R}(P_R(\lam),\Delta_R(\mu))   @>J>> \Hom_{\OC_R}(P_R(\lam),\nabla_R(\mu))\\
   @VVt V @VVL V\\
   \Hom_{\OC_R}(\Delta_R(\mu'),K_R(\lam')) @>E >> \left(\Hom_{\OC_R}(K_R(\lam'),\nabla_R(\mu'))\right)^{*}
	\end{CD}
\end{eqnarray}
commutative. Here $J=J_R^{\mu}(P_R(\lam))$ and $E=E_R^{\mu'}(K_R(\lam'))$ denote the inclusions defined above and $t=t_R$ denotes the isomorphism induced by the tilting functor.
\end{thm}

\begin{proof}
If $\lam$ is not contained in the equivalence class of $\mu$ under $\sim_R = \sim_\C$ all $\Hom$-spaces occurring in the diagram are 0 by block decomposition and the assertion of the proposition is true.\\
So let us assume $\lam$ to be in the equivalence class of $\mu$. It is easy to see that $J$ and $E$ commute with base change and we have also verified this property for $t$ in Lemma \ref{tilt} already. Let $\p \subset R$ be a prime ideal of height 1. We abbreviate $\Hom=\Hom_{\OC_{R_\p}}, P=P_R(\lam), K=K_R(\lam'), \Delta=\Delta_{R_\p}(\mu), \nabla=\nabla_{R_\p}(\mu), \Delta'=\Delta_{R_\p}(\mu')$ and $\nabla'=\nabla_{R_\p}(\mu')$. After applying $\cdot \otimes_R R_\p$ to our diagram and the base change isomorphisms of Theorem \ref{Fie2} and Corollary \ref{Soe2} we get a diagram of $R_\p$-modules

\begin{eqnarray}\label{eqn:Dia}
	\begin{CD}
   \Hom(P\otimes_{R}R_{\p},\Delta)   @>J>> \Hom(P\otimes_{R}R_{\p},\nabla)\\
   @VVt V\\
   \Hom(\Delta',K\otimes_{R}R_{\p}) @>E >> \left(\Hom(K\otimes_{R}R_{\p},\nabla')\right)^{*_{R_\p}}
	\end{CD}
\end{eqnarray}
where we omit the index $R_\p$ for $t, J$ and $E$.\\
We want to show that we get an isomorphism $L_{R_\p}$ of $R_\p$-modules as the missing right vertical of the upper diagram to make it commutative for every prime ideal $\p \subset R$ of height 1. By block decomposition we get
$$P_R(\lam) \otimes_R R_\p \cong \bigoplus\limits_{i=1}^{n} P_{R_\p}(\lam_i)$$
for certain indecomposable projective objects $P_{R_\p}(\lam_i) \in \OC_{R_\p}$ and $\lam_i \in \overline{\lam}$ where $\overline{\lam}$ denotes the equivalence class of $\lam$ under $\sim_R=\sim_\C$. Since $t$ is fully faithful and respects base change we also get a decomposition
$$K_R(\lam') \otimes_R R_\p \cong t_{R_\p}(P_R(\lam) \otimes_R R_\p) \cong \bigoplus\limits_{i=1}^{n} t_{R_\p}(P_{R_\p}(\lam_i))$$
It is easy to see that $J$, $E$ and $t$ respect these decompositions. Hence, we get in formulas

\begin{eqnarray*}
	\begin{array}{ccc}
  J(\Hom_{\OC_{R_{\mathfrak{p}}}}(P_{R_{\mathfrak{p}}}(\lam_i),\Delta_{R_{\mathfrak{p}}}(\mu)))&\subset & \Hom_{\OC_{R_{\mathfrak{p}}}}(P_{R_{\mathfrak{p}}}(\lam_i),\nabla_{R_{\mathfrak{p}}}(\mu))\\
  \\
  t(\Hom_{\OC_{R_{\mathfrak{p}}}}(P_{R_{\mathfrak{p}}}(\lam_i),\Delta_{R_{\mathfrak{p}}}(\mu)))&=& \Hom_{\OC_{R_{\mathfrak{p}}}}(t(\Delta_{R_{\mathfrak{p}}}(\mu)),t(P_{R_{\mathfrak{p}}}(\lam_i)))\\
  \\
  E\left(\Hom_{\OC_{R_{\mathfrak{p}}}}(t(\Delta_{R_{\mathfrak{p}}}(\mu)),t(P_{R_{\mathfrak{p}}}(\lam_i)))\right)&\subset& \left(\Hom_{\OC_{R_{\mathfrak{p}}}}(t(P_{R_{\mathfrak{p}}}(\lam_i)),\nabla_{R_{\mathfrak{p}}}(\mu'))\right)^{*_{R_\p}}
	\end{array}
\end{eqnarray*}
where we omit the index $R_\p$ of our maps.\\
Let $Q=\text{Quot}(S)$ be the quotient field of $S$. Since all deformed Verma modules over $Q$ are simple, we get $\Delta_Q(\mu) \cong \nabla_Q(\mu)$ and $\Delta_Q(\mu') \cong \nabla_Q(\mu')$, respectively. Hence, after applying $\cdot \otimes_R Q$ to our diagram we get an isomorphism $L_Q$ which makes the following diagram commutative

{\small
\begin{eqnarray*}
	\begin{CD}
   \Hom_{\OC_Q}(P_R(\lam)\otimes_RQ,\Delta_Q(\mu))   @>J_Q>> \Hom_{\OC_Q}(P_R(\lam)\otimes_RQ,\nabla_Q(\mu))\\
   @VVt_Q V  @VVL_QV\\
   \Hom_{\OC_Q}(t_Q(\Delta_Q(\mu)),t_Q(P_R(\lam)\otimes_RQ)) @>E_Q >> \left(\Hom_{\OC_Q}(t_Q(P_R(\lam)\otimes_RQ),\nabla_Q(\mu'))\right)^{*_Q}
	\end{CD}
\end{eqnarray*}}
Let $\p \subset R$ be a prime ideal of height one. So $\p = R\gamma$ for an irreducible element $\gamma \in R$. If $\gamma \notin \C \check{\alpha}$ for all $\alpha \in \mathcal{R}^{+}$ or $\gamma \in \C\check{\alpha}$ for one $\alpha \in \mathcal{R}$ but with $\langle \mu + \rho, \check{\alpha}\rangle_\C \notin \Z$, the block of $\mu$ in $\OC_{R_\p}$ is generic by Lemma \ref{Fie4} and we conclude 
$P_{R_\mathfrak{p}}(\mu) \cong \Delta_{R_\mathfrak{p}}(\mu) \cong \nabla_{R_\mathfrak{p}}(\mu)$. But in this case all maps $J_{R_\p}$, $E_{R_\p}$ and $t_{R_\p}$ are isomorphisms and so we get the claimed ${R_\p}$-isomorphism $L_{R_\p}$ of our diagram (\ref{eqn:Dia}).
So we can assume $\p = R\check{\alpha}$ for some $\alpha \in \mathcal{R}^{+}$ and $\langle \mu + \rho, \check{\alpha}\rangle_\C \in \Z$.\\
Now we want to show 

{\small
\begin{eqnarray*}
	\begin{array}{rcl}
 L_Q(\Hom_{\OC_{R_\mathfrak{p}}}(P_{R}(\lam)\otimes_{R}{R_\mathfrak{p}},\nabla_{R_\mathfrak{p}}(\mu)))&=&\left(\Hom_{\OC_{R_\mathfrak{p}}}(t(P_{R}(\lam))\otimes_{R}{R_\mathfrak{p}},\nabla_{R_\mathfrak{p}}(\mu'))\right)^{*_{R_\mathfrak{p}}}\\
 \\
 
 &\subset& \left(\Hom_{\OC_Q}(t(P_{R}(\lam))\otimes_{R}Q,\nabla_Q(\mu'))\right)^{*_Q}
	\end{array}
\end{eqnarray*}}
Since all our maps respect the decomposition of $P_R(\lam) \otimes_R R_\p$, we only have to prove this for the indecomposable summands of our decomposition. In formulas we want to show

{\footnotesize \begin{eqnarray}\label{eqn:ZuZe}
	\begin{array}{rcl}
 L_Q(\Hom_{\OC_{R_\mathfrak{p}}}(P_{R_\mathfrak{p}}(\lam_i),\nabla_{R_\mathfrak{p}}(\mu)))&=&\left(\Hom_{\OC_{R_\mathfrak{p}}}(t(P_{R_\mathfrak{p}}(\lam_i)),\nabla_{R_\mathfrak{p}}(\mu'))\right)^{*_{R_\mathfrak{p}}}
	\end{array}
\end{eqnarray}}
for all $i \in \{1,...,n\}$. Here, the right side is again meant to be the $R_\p$-lattice in the $Q$-vector space $\Hom_{\OC_Q}(t(P_{R_\p}(\lam_i))\otimes_{R_\p} Q, \nabla_{Q}(\mu'))^{*_Q}$.\\
Both $\Hom$-spaces in (\ref{eqn:ZuZe}) are free and of the same rank over $R_\p$. From the description of the projective covers in the generic and subgeneric case it follows
\begin{eqnarray*}
	\begin{array}{rcl}
\text{rk}_{R_\mathfrak{p}}(\Hom_{\OC_{R_\mathfrak{p}}}(P_{R_\mathfrak{p}}(\lam_i),\nabla_{R_\mathfrak{p}}(\mu)))&=&\text{dim}_{\K}(\Hom_{\OC_{\K}}(P_{\K}(\lam_i),\nabla_{\K}(\mu)))\\
\\
 &=& (P_{\K}(\lam_i):\Delta_{\K}(\mu)) \leq 1
	\end{array}
\end{eqnarray*}
where $\K$ denotes the residue field of $R_\p$.\\
Now we will proceed in several steps to prove (\ref{eqn:ZuZe}). From the choice of $\p$ it follows that the equivalence class of $\mu$ under $\sim_{R_\p}$ equals $\{\mu, s_\alpha \cdot \mu\}$ where again $s_\alpha$ is the reflection corresponding to the root $\alpha \in \mathcal{R}$ acting by the dot-action on $\h^{*}$. If $\mu = s_\alpha \cdot \mu$, we are in the generic case and again our maps are isomorphisms which proves the claim for this case. So let us further assume $\mu \neq s_\alpha \cdot \mu$.

\begin{stepone} 
Let $\lam_i \notin \{\mu, s_\alpha \cdot \mu\}$.\\
In this case both $\Hom$-spaces in (\ref{eqn:ZuZe}) are zero by block decomposition and the claim is true.
\end{stepone}

\begin{steptwo} 
Let $\lam_i =s_\alpha \cdot \mu > \mu$.\\
In this case we get $P_{R_{\mathfrak{p}}}(\lam_i) \cong \Delta_{R_{\mathfrak{p}}}(\lam_i)\neq \Delta_{R_{\mathfrak{p}}}(\mu) \cong \nabla_{R_{\mathfrak{p}}}(\mu)$ and therefore
\begin{eqnarray*}
	\begin{array}{rcl}
	\Hom_{\OC_{R_\mathfrak{p}}}(P_{R_\mathfrak{p}}(\lam_i),\Delta_{R_\mathfrak{p}}(\mu))&=&\Hom_{\OC_{R_\mathfrak{p}}}(P_{R_\mathfrak{p}}(\lam_i),\nabla_{R_\mathfrak{p}}(\mu))\\
	\\
	&=& \Hom_{\OC_{R_\mathfrak{p}}}(t(\Delta_{R_\mathfrak{p}}(\mu)),t(P_{R_\mathfrak{p}}(\lam_i)))\\
	\\
	&=&\left(\Hom_{\OC_{R_\mathfrak{p}}}(t(P_{R_\mathfrak{p}}(\lam_i)),\nabla_{R_\mathfrak{p}}(\mu'))\right)^{*}=0
	\end{array}
\end{eqnarray*}
and we are done.

\end{steptwo}

\begin{stepthree} 
Let $\lam_i = \mu< s_\alpha \cdot \mu$.\\
Then we have $\Delta_{R_\mathfrak{p}}(\mu)=\nabla_{R_\mathfrak{p}}(\mu)$ and also $t(P_{R_\mathfrak{p}}(\lam_i))\cong K_{R_\mathfrak{p}}(-2\rho-\lam_i)\cong P_{R_\mathfrak{p}}(-2\rho-s_\alpha \cdot \mu)$. With these assumptions on $\lam_i$ and $\mu$, we get an injection 
$$\varphi :t(\Delta_{\K}(\mu)) \hookrightarrow t(P_{\K}(\lam_i))$$
Since $-2 \rho - s_\alpha \cdot \mu < -2 \rho - \mu$, the map $\varphi$ is an embedding from a Verma module with dominant highest weight to the projective cover of the simple module with antidominant highest weight. These modules are objects of a subgeneric block of $\OC_\K$ and this projective cover is self-dual. So dualizing $\varphi$ leads to a surjection $d \varphi: P_\K(-2\rho -s_\alpha \cdot \mu) \twoheadrightarrow \nabla_\K(-2\rho - \mu)$. But the composition $d\varphi \circ \varphi$ is nonzero. Thus, we get $E_\K(\varphi)\neq 0$. Comparing dimensions over $\K$, we get an isomorphism
\begin{eqnarray*}
	\begin{array}{c}
 E_\K :\Hom_{\OC_{\K}}(t(\Delta_{\K}(\mu)),t(P_{\K}(\lam_i)))\stackrel{\sim}{\rightarrow}	\left(\Hom_{\OC_{\K}}(t(P_{\K}(\lam_i)),\nabla_{\K}(\mu'))\right)^{*}
		\end{array}
\end{eqnarray*}
By Nakayama's lemma, we conclude that $E_{R_\p}$ is surjective and therefore an isomorphism.\\
We also have a surjective map $P_\K(\lam_i) \twoheadrightarrow \Delta_\K(\mu)$ and with similar arguments we get an isomorphism of vector spaces
\begin{eqnarray*}
	\begin{array}{c}
\Hom_{\OC_{\K}}(P_{\K}(\lam_i),\Delta_{\K}(\mu))\stackrel{\sim}{\rightarrow}	\Hom_{\OC_{\K}}(P_{\K}(\lam_i),\nabla_{\K}(\mu))
		\end{array}
\end{eqnarray*}
and Nakayama's lemma finishes this case.
\end{stepthree}

\begin{stepfour} 
Let $\lam_i = \mu> s_\alpha \cdot \mu$.\\
Then we have  $P_{R_\mathfrak{p}}(\lam_i)= \Delta_{R_\mathfrak{p}}(\lam_i)=\Delta_{R_\mathfrak{p}}(\mu)$ and over the residue field $\K$ we get $E_{\K}(\text{id})\neq 0$ and $J_{\K}(\text{id})\neq 0$. And again by Nakayama's lemma we get two isomorpisms

\begin{eqnarray*}
	\begin{array}{clcl}
J:&\Hom_{\OC_{R_\mathfrak{p}}}(P_{R_\mathfrak{p}}(\lam_i),\Delta_{R_\mathfrak{p}}(\mu))&\stackrel{\sim}{\rightarrow}&	\Hom_{\OC_{R_\mathfrak{p}}}(P_{R_\mathfrak{p}}(\lam_i),\nabla_{R_{\mathfrak{p}}}(\mu))\\
\\
E:& \Hom_{\OC_{R_\mathfrak{p}}}(t(\Delta_{R_\mathfrak{p}}(\mu)),t(P_{R_\mathfrak{p}}(\lam_i)))&\stackrel{\sim}{\rightarrow}&	\left(\Hom_{\OC_{R_\mathfrak{p}}}(t(P_{R_\mathfrak{p}}(\lam_i)),\nabla_{R_{\mathfrak{p}}}(\mu'))\right)^{*_{R_\p}}
		\end{array}
\end{eqnarray*}
and the claim is true.
\end{stepfour}

\begin{stepfive} 
Let $\lam_i =  s_\alpha \cdot \mu < \mu$ and set $\lam_i'=-2\rho-\lam_i$.\\
After base change with the residue field $\K$ we get
$$J_{\K}:\Hom_{\OC_{\K}}(P_{\K}(\lam_i),\Delta_{\K}(\mu))\longrightarrow	\Hom_{\OC_{\K}}(P_{\K}(\lam_i),\nabla_{\K}(\mu))$$
On the left side we have the generator $P_{\K}(\lam_i) \twoheadrightarrow \Delta_{\K}(\lam_i) \hookrightarrow \Delta_{\K}(\mu)$ which becomes
$$P_{\K}(\lam_i) \twoheadrightarrow \Delta_{\K}(\lam_i) \hookrightarrow \Delta_{\K}(\mu)\twoheadrightarrow L_{\K}(\mu) \hookrightarrow \nabla_{\K}(\mu) $$
on the right side. But since this map is zero we get
$$J\left(\Hom_{\OC_{R_\mathfrak{p}}}(P_{R_\mathfrak{p}}(\lam_i),\Delta_{R_\mathfrak{p}}(\mu))\right)\subset \check{\alpha} \cdot \Hom_{\OC_{R_\mathfrak{p}}}(P_{R_\mathfrak{p}}(\lam_i),\nabla_{R_{\mathfrak{p}}}(\mu))$$
Now we can lift the generator
$j':P_{\K}(\lam_i) \twoheadrightarrow \Delta_{\K}(\lam_i) \hookrightarrow \Delta_{\K}(\mu)$ to a generator\\ $j:P_{R_\mathfrak{p}}(\lam_i) \rightarrow \Delta_{R_\mathfrak{p}}(\mu)$ of the $R_\mathfrak{p}$-module $\Hom_{\OC_{R_{\mathfrak{p}}}}(P_{R_{\mathfrak{p}}}(\lam_i),\Delta_{R_{\mathfrak{p}}}(\mu))$ via base change. Let also $l\in \Hom_{\OC_{R_{\mathfrak{p}}}}(P_{R_{\mathfrak{p}}}(\lam_i),\nabla_{R_{\mathfrak{p}}}(\mu))$ be a generator. Then we conclude $J_{R_{\mathfrak{p}}}(j)=cl$ for one $c\in {R_{\mathfrak{p}}}$ with $c \neq 0$. Multiplying with an appropriate invertible element of $R_\p$ we can assume $c=(\check{\alpha})^{m}$ for $m \in \N$.\\
On $\Delta_{R_{\mathfrak{p}}}(\mu)$ we have the filtration which comes from the preimages of $(\check{\alpha})^{i} \nabla_{R_{\mathfrak{p}}}(\mu)$ under the inclusion $\Delta_{R_{\mathfrak{p}}}(\mu) \hookrightarrow \nabla_{R_{\mathfrak{p}}}(\mu)$. Denote by $\Delta_{\K}(\mu)_i$ the image of the i-th filtration layer under the surjection $\Delta_{R_{\mathfrak{p}}}(\mu) \twoheadrightarrow \Delta_\K(\mu)$ which is the i-th layer of the Jantzen filtration on $\Delta_\K(\mu)$. By the Jantzen sum formula we get
$$\sum_{i>0} \text{ch} \Delta_{\K}(\mu)_i = \sum_{\substack{\alpha \in \mathcal{R}^{+}, n\in \N\\  \langle \mu + \rho +\tau, \check{\alpha} \rangle_{\K}=n}}\text{ch} \Delta_{\K}(\mu -n\alpha)=\text{ch}\Delta_{\K}(s_\alpha \cdot \mu)$$
Since $\Delta_{\K}(s_\alpha \cdot \mu)$ is simple, we conclude $\Delta_{\K}(\mu)_1\cong \Delta_{\K}(s_\alpha \cdot \mu)$ and $\Delta_{\K}(\mu)_2=0$.\\
If $J_{R_\p}(j)\in (\check{\alpha})^{2}\cdot \Hom_{\OC_{R_\mathfrak{p}}}(P_{R_\mathfrak{p}}(\lam_i),\nabla_{R_{\mathfrak{p}}}(\mu))$ we would get $j'(P_{\K}(\lam_i))\subset \Delta_{\K}(\mu)_2=0$, i.e. $j'=0$, which is a contradiction to the choice of $j'$. We conclude $J_{R_\p}(j)=\check{\alpha}l$.\\

Now we will observe the lower horizontal of the diagram. We have $t_{R_\p}(P_{R_\p}(\lam_i))=K_{R_\p}(\lam_i')$ and by our assumptions we get $s_\alpha \cdot (\lam_i')=\mu' <\lam_i'$. But since $\mu'$ is minimal in its equivalence class under $\sim_{R_\p}$, we get an isomorphism $P_{\K}(\mu')\cong K_{\K}(\lam_i')$ and then an inclusion $\gamma':\Delta_{\K}(\mu')\hookrightarrow K_{\K}(\lam_i')$. Lifting this map with base change \ref{Soe2}, we get a basis $\gamma$ of the $R_\p$-module $\Hom_{\OC_{R_\p}}(\Delta_{R_\p}(\mu'),K_{R_\p}(\lam_i'))$. But every map $\Delta_{\K}(\mu') \rightarrow \nabla_{\K}(\mu')$ which factors through $K_{\K}(\lam_i')$ is zero already. So we conclude that
\begin{displaymath}
\begin{array}{c}
E_{\K}: \Hom_{\OC_{\K}}(\Delta_{\K}(\mu'),K_{\K}(\lam_i')) \rightarrow
 \left(\Hom_{\OC_{\K}}(K_{\K}(\lam_i'),\nabla_{\K}(\mu')) \right)^{*}
\end{array}
\end{displaymath}
is the zero map. If we now choose a basis $\delta$ of of the free one dimensional $R_\p$-module $\left(\Hom_{\OC_{R_\p}}(K_{R_\p}(\lam_i'),\nabla_{R_\p}(\mu')) \right)^{*}$ and eventually multiply it with an appropriate invertible element of $R_\p$, we conclude $E_{R_\p}(\gamma)=(\check{\alpha})^n \delta$ for $n>0$.\\
By Proposition \ref{Fie5} we get

\begin{displaymath}
\begin{array}{rcl}
\End_{\OC_{R_\p}}(K_{R_\p}(\lam_i'))&\cong&\End_{\OC_{R_\p}}(P_{R_\p}(\mu'))\\
 &\cong& \left \{(x,y) \in R_{\mathfrak{p}}\oplus R_{\mathfrak{p}} \middle| x \equiv y \text{ mod } \check{\alpha}\right\}
\end{array}
\end{displaymath}
where for $\varphi \in \End_{\OC_{R_\p}}(K_{R_\p}(\lam_i'))$ $x$ and $y$ are given by the induced maps on the short exact sequences

\begin{eqnarray*}
\begin{CD}
   0   @>>> \Delta_{R_\p}(\lam'_i) @>>> K_{R_\p}(\lam'_i) @>>> \Delta_{R_\p}(\mu') @>>> 0\\
   @VVV @Vx \cdot VV @VV\varphi V @VVy \cdot V @VVV\\
   0   @>>> \Delta_{R_\p}(\lam'_i) @>>> K_{R_\p}(\lam'_i) @>>> \Delta_{R_\p}(\mu') @>>> 0 
\end{CD}\end{eqnarray*}

If for $\varphi$ we choose the homomorphism that corresponds to the tupel $(x,y)=(0,\check{\alpha})$, the map on the cokernels of the diagram will factor through the middle and we get a map
$$\Delta_{R_\p}(\mu') \rightarrow K_{R_\p}(\lam'_i)\rightarrow  \Delta_{R_\p}(\mu')=\nabla_{R_\p}(\mu')$$
and the image of this map is $\check{\alpha}\nabla_{R_\p}(\mu')$. So we get $n=1$. This proves the last case.\\

For now, we showed that for every prime ideal $\p \subset R$ of height 1 we find an isomorphism $L_{R_\p}$ that makes the diagram (\ref{eqn:Dia}) commutative. Since $E$, $J$ and $t$ respect base change we get an isomorphism which we will also call $L=L_{R_\p}$ that makes 
{\small
\begin{eqnarray*}
\begin{CD}
   \Hom_{\OC_{R}}(P_{R}(\lam),\Delta_{R}(\mu))\otimes_{R}R_\p @>J>> \Hom_{\OC_{R}}(P_{R}(\lam),\nabla_{R}(\mu))\otimes_{R}R_\p \\
   @VtVV @VVLV \\
   \Hom_{\OC_{R}}(t(\Delta_{R}(\mu)),t(P_{R}(\lam)))\otimes_{R}R_\p @>E>> \left(\Hom_{\OC_{R}}(t(P_{R}(\lam)),\nabla_{R}(\mu'))\right)^{*} \otimes_{R}R_\p
\end{CD}\end{eqnarray*}}
commutative.\\
But then we also get an isomorphism $L_{R}$ which is the restriction of the $Q$-linear map $L_Q$ to $\Hom_{\OC_{R}}(P_{R}(\lam),\nabla_{R}(\mu))$ and makes the diagram
{\footnotesize
\begin{eqnarray*}
\begin{CD}
   \bigcap\limits_{\p \in \mathfrak{P}}\Hom_{\OC_{R}}(P_{R}(\lam),\Delta_{R}(\mu))\otimes_{R}R_\p @>J>> \bigcap\limits_{\p \in \mathfrak{P}}\Hom_{\OC_{R}}(P_{R}(\lam),\nabla_{R}(\mu))\otimes_{R}R_\p \\
   @VtVV @VVLV \\
   \bigcap\limits_{\p \in \mathfrak{P}}\Hom_{\OC_{R}}(t(\Delta_{R}(\mu)),t(P_{R}(\lam)))\otimes_{R}R_\p @>E>> \bigcap\limits_{\p \in \mathfrak{P}}\left(\Hom_{\OC_{R}}(t(P_{R}(\lam)),\nabla_{R}(\mu'))\right)^{*} \otimes_{R}R_\p
\end{CD}
\end{eqnarray*}}
commutative, where $\mathfrak{P}$ denotes the set of all prime ideals of $R$ of height one. Since all $\Hom$-spaces of this diagram are finitely generated free $R$-modules, it equals diagram (\ref{eqn:DiaS}).
\end{stepfive}
\end{proof}

\begin{cor}
Let $P \in \OC$ be a projective object, $\mu \in \h^{*}$ and let $t=t_\C$ be the tilting functor. Then $K=t(P)$ is a tilting object in $\OC$ and the isomorphism\\ $t:\Hom_\g(P,\Delta(\mu)) \stackrel{\sim}{\rightarrow} \Hom_\g(\Delta(-2\rho - \mu),K)$ induced by the tilting functor identifies the filtration induced by the Jantzen filtration with the Andersen filtration.
\end{cor}

\begin{proof}
Consider the restriction map $S=\C[\h^{*}] \rightarrow \C \llbracket v \rrbracket$ to a formal neighbourhood of the origin in the line $\C\rho$. Since this map factors through $R$, we get a homomorphism $R \rightarrow \C \llbracket v \rrbracket$ of $S$-algebras. Since the maps from the proposition above commute with base change $\cdot \otimes_R \C \llbracket v \rrbracket$, we get a commuting diagram

\begin{eqnarray*}
	\begin{CD}
   \Hom_{\OC_{\C \llbracket v \rrbracket}}(P_{\C \llbracket v \rrbracket},\Delta_{\C \llbracket v \rrbracket}(\mu))   @>J>> \Hom_{\OC_{\C \llbracket v \rrbracket}}(P_{\C \llbracket v \rrbracket},\nabla_{\C \llbracket v \rrbracket}(\mu))\\
   @VVt V @VVL V\\
   \Hom_{\OC_{\C \llbracket v \rrbracket}}(\Delta_{\C \llbracket v \rrbracket}(\mu'),K_{\C \llbracket v \rrbracket}) @>E >> \left(\Hom_{\OC_{\C \llbracket v \rrbracket}}(K_{\C \llbracket v \rrbracket},\nabla_{\C \llbracket v \rrbracket}(\mu'))\right)^{*}
	\end{CD}
\end{eqnarray*}
where $P_{\C \llbracket v \rrbracket}$ (resp. $K_{\C \llbracket v \rrbracket}$) denotes the unique (up to isomorphism) deformed projective (resp. tilting) module in $\OC_{\C \llbracket v \rrbracket}$ that maps to $P$ (resp. $K$) under $\cdot \otimes_{\C \llbracket v \rrbracket} \C$, $\mu' =-2\rho- \mu$ and $(\cdot)^{*}$ means the $\C \llbracket v \rrbracket$-dual. Since $L(v^{i}\Hom_{\OC_{\C \llbracket v \rrbracket}}(P_{\C \llbracket v \rrbracket},\nabla_{\C \llbracket v \rrbracket}(\mu)))=v^{i}(\Hom_{\OC_{\C \llbracket v \rrbracket}}(K_{\C \llbracket v \rrbracket},\nabla_{\C \llbracket v \rrbracket}(\mu')))^{*}$, the preimages of these submodules under $J$ resp. $E$ are identified by $t$. Applying $\cdot \otimes_{\C \llbracket v \rrbracket} \C$ to the left vertical of the diagram yields the claim.
\end{proof}
The last corollary shows the connection between the Jantzen and Andersen filtrations. Taking a projective generator of the block containing the Verma module $\Delta(\mu)$, the above filtration on the $\C$-vector space $\Hom_\g(P,\Delta(\mu))$ might carry enough information to get back the Jantzen filtration on $\Delta(\mu)$. If the map $t$ between the two Hom-spaces induced by the tilting functor also respects the gradings coming from the graded version of category $\OC$, the results of \cite{I} about the Andersen filtration could give an alternative proof of the Jantzen conjecture about the semisimplicity of the subquotients of the Jantzen filtration which was proved in \cite{A}.


\bibliographystyle{amsplain}

\end{document}